\documentclass[reqno,centertags, 12pt]{amsart}
\usepackage{amscd}
\usepackage{latexsym}
\usepackage{amsmath}
\usepackage{amsthm}
\usepackage{amsfonts}
\usepackage{amssymb}
\numberwithin{equation}{section}

\newcommand{\vp}{\varphi}
\newcommand{\T}{\partial\mathbb{D}}

\newcommand{\ds}{\displaystyle}
\newcommand{\ol}{\overline}

\newcommand{\be}{\begin{equation}}
\newcommand{\ee}{\end{equation}}
\newcommand{\ba}{\begin{array}}
\newcommand{\ea}{\end{array}}
\newcommand{\C}{\mathbb{C}}

\newcommand{\la}{\langle}
\newcommand{\ra}{\rangle}

\newtheorem{lemma}{Lemma}[section]
\newtheorem{theorem}{Theorem}[section]

\begin{document}

\title[A Formula for Inserting Point Masses]
{A Formula for Inserting Point Masses}
\author[M.-W. L. Wong]{Manwah Lilian Wong}
\thanks{$^*$ Mathematics 253-37, California Institute of Technology, Pasadena, CA 91125.
E-mail: wongmw@caltech.edu. Supported by the Croucher Foundation, Hong Kong}
\date{October 15th, 2007}
\keywords{point masses, decay of Verblunsky coefficients}
\subjclass[2000]{42C05, 30E10, 05E35}

\begin{abstract} Let $d\mu$ be a probability measure on the unit circle and $d\nu$ be the measure formed by adding a pure point to $d\mu$.  We give a formula for the Verblunsky coefficients of $d\nu$, based on a result of Simon.
\end{abstract}
\maketitle

\section{Introduction}

Suppose we have a probability measure $d\mu$ on the unit circle $\T=\{z \in \C: |z|=1\}$. We define the inner product associated with $d\mu$ and the norm on $L^2(\T, d\mu)$ respectively by
\begin{align}
\left\langle f ,g \right\rangle & = \ds \int_{\T} \ol{f(e^{i \theta})} g(e^{i \theta}) d\mu(\theta) \\
\|f\|_{d\mu} & = \left( \ds \int_{\T} |f(e^{i \theta})|^2 d\mu(\theta) \right)^{1/2}
\end{align}

The family of monic orthogonal polynomials associated with the measure $d\mu$ is denoted as $(\Phi_n(z, d\mu))_{n=0}^{\infty}$, while the normalized family is denoted as $(\vp_n(z, d\mu))_{n=0}^{\infty}$.

Let $\Phi_n^*(z)=z^n \ol{\Phi_n(1/\ol{z})}$ and $\vp_n^*(z)=\Phi_n^*(z)/\|\Phi_n\|$ be the reversed polynomials. Orthogonal polynomials obey the Szeg\H o recursion relation
 \be
\Phi_{n+1}(z) = z \Phi_n (z) - \ol{\alpha_n} \Phi_n^*(z)
\label{eq01}
\ee
$\alpha_n$ is called the $n^{th}$ Verblunsky coefficient. It is well known that there is a one-to-one correspondence between $d\mu$ and $(\alpha_j(d\mu))_{j=0}^{\infty}$ and that the Verblunsky coefficients carry much information about the family of orthogonal polynomials. For example, 
\be
\| \Phi_n\|^2 = \ds \prod_{j=0}^{n-1} (1-|\alpha_j|^2)
\label{normphi}
\ee

For a comprehensive introduction to the theory of orthogonal polynomials on the unit circle, the reader should refer to \cite{simon1, simon2}, or the classic reference \cite{szego}.
\medskip

The result that we would like to present is the following
\begin{theorem} Suppose $d\mu$ is a probability measure on the unit circle and $0<\gamma<1$. Let $d\nu$ be the probability measure formed by adding a point mass $\zeta = e^{i \omega} \in \T$ to $d\mu$ in the following manner
\be
d\nu = (1-\gamma) d\mu + \gamma \delta_{\omega}
\label{dnu}
\ee
Then the Verblunsky coefficients of $d\nu$ are given by
\be
\alpha_n(d\nu) = \alpha_n + \ds \frac{(1-|\alpha_n|^2)^{1/2}}{(1-\gamma) \gamma^{-1}+ K_n(\zeta)} \ol{\vp_{n+1}(\zeta)} \vp_n^*(\zeta)
\label{myformula}
\ee where
\be
K_n(\zeta) = \ds \sum_{j=0}^{n} |\vp_j(\zeta)|^2
\ee and all objects without the label $(d\nu)$ are associated with the measure $d\mu$.
\label{theorem0}
\end{theorem}

The proof is based on a result obtained by Simon in the proof of Theorem 10.13.7 in \cite{simon2} (See Theorem \ref{barrytheorem} below).

In fact, the following formula had been found by Geronimus \cite{geronimus}
\be
\Phi_{n}(z, d\nu) = \Phi_n(z) - \ds \frac{\Phi_n(\zeta) K_{n-1}(z, \zeta)}{(1-\gamma) \gamma^{-1} + K_{n-1}(\zeta, \zeta)}
\label{geronimus}
\ee

The formula for the real case was rediscovered by Nevai \cite{nevai}, Later, the same formula for the unit circle case was rediscovered by Cachafeiro-Marcellan \cite{cm}. Unaware of Geronimus' result and the fact that Nevai's result also applies to the unit circle, Simon reconsidered this problem and proved formula (\ref{mainformula}) independently using a totally different method.

For applications of formula (\ref{myformula}), the reader may refer to \cite{wong1} and \cite{wong2}.

\medskip

\section{The Proof}\label{proof0}
First, we will prove a few lemmas.
\begin{lemma} Let $\beta_{jk}=\la \Phi_j(d\mu), \Phi_k (d\mu) \ra_{d\nu}$. Then
\be
\Phi_n(d\nu)(z)=\ds \frac{1}{D^{(n-1)}} \left|
\begin{matrix}
\beta_{00} & \beta_{0 \,1} & \dots & \beta_{0 \,n} \\ 
\vdots & & & \vdots \\
\beta_{{n-1}\, 0} & \beta_{{n-1}\, 1} & \dots & \beta_{{n-1} \, n}\\
\Phi_0 (d\mu) & \dots &  \dots & \Phi_n(d\mu)
\end{matrix} \right|
\label{rhs}
\ee where
\be
D^{(n-1)} =  \left|
\begin{matrix}
\beta_{0\,0} & \beta_{0\,1} & \dots & \beta_{0 \,n-1} \\ 
\vdots & & & \vdots \\
\beta_{{n-1}\, 0} & \beta_{{n-1}\, 1} & \dots & \beta_{{n-1}\, {n-1}}\\
\end{matrix} \right|
\ee
\label{lemma1}
\end{lemma}

\begin{proof}
Let $\tilde{\Phi}_n(d\nu)$ be the right hand side of (\ref{rhs}). We observe that the inner product $\langle \Phi_j(d\mu), \tilde{\Phi}_n(d\nu)\rangle_{d\nu}$ is zero for $j=0,1, \dots, n-1$ as the last row and the $j^{th}$ row of the determinant are the same. By expanding in minors, we see that the leading coefficient of $\tilde{\Phi}_n(d\nu)$ in (\ref{rhs}) is one. In other words, $\tilde{\Phi}_n(d\nu)$ is an $n^{th}$ degree monic polynomial which is orthogonal to $1, z, \dots, z^{n-1}$ with respect to $\langle \,,\, \rangle_{d\nu}$, hence $\tilde{\Phi}_n(d\nu)$ equals $\Phi_n(d\nu)$.
\end{proof}

\begin{lemma} Let $C$ be the following $(n+1) \times (n+1)$ matrix
\be
\begin{pmatrix} A& v \\ w & \beta \end{pmatrix} 
\ee where $A$ is an $n \times n$ matrix, $\beta$ is in $\mathbb{C}$, $v$ is the column vector $(v_0, v_1, \dots, v_{n-1})^{T}$ and $w$ is the row vector $(w_0, w_1, \dots, w_{n-1})$. If $det(A) \not =0$, we have
\be
\det(C) = \det(A) \left( \beta  - \ds \sum_{0 \leq j,k \leq n-1} w_k v_j (A^{-1})_{j k}  \right)
\ee
\label{lemma2}
\end{lemma}

\begin{proof}
We expand in minors, starting from the bottom row to get
\be
\det(C) = \beta det(A)  + \ds \sum_{0 \leq j,k \leq n-1} w_k v_j (-1)^{j+k+1} \det(\tilde{A}_{jk})
\ee where $\tilde{A}_{jk}$ is the matrix $A$ with the $j^{th}$ row and $k^{th}$ column removed.

By Cramer's rule, since $\det(A) \not = 0$,
\be
\tilde{A}_{jk}= (-1)^{j+k} \det(A) (A^{-1})_{j k}
\ee proving Lemma \ref{lemma2}.
\end{proof} 

Next, we are going to prove the following formula by Simon \cite{simon2}:
\begin{theorem} The Verblunsky coefficient of $d\nu$ (as defined in (\ref{dnu})) is given by
\begin{equation}
\alpha_{n}(d\nu)=\alpha_{n}-q_{n}^{-1} \gamma \ol{\vp_{n+1}(\zeta)} \left( \ds \sum_{j=0}^{n} \alpha_{j-1} \frac{\|\Phi_{n+1}\|}{\|\Phi_j\|}\vp_j(\zeta) \right)
\label{mainformula}
\end{equation} where
\begin{align}
K_{n}(\zeta) & = \ds \sum_{j=0}^{n} |\vp_j(\zeta)|^2 \label{kndef}\\
q_{n} & = (1-\gamma) + \gamma K_{n}(\zeta) \\
\alpha_{-1} & = -1
\end{align} and all objects without the label $(d\nu)$ are associated with the measure $d\mu$.
\label{barrytheorem}
\end{theorem}

\begin{proof}  Since $\alpha_{n-1} (d\nu) =  - \ol{\Phi_n(0, d\nu)}$ and $\ol{\beta_{jk}}=\beta_{kj}$, by Lemma \ref{lemma1},
\be
\alpha_{n-1}(d\nu) 
%=\ds \frac{1}{\ol{D^{(n-1)}}} \left|
%\begin{matrix}
%\ol{\beta_{00}} & \ol{\beta_{0 1}} & \dots & \ol{\beta_{0 n}} \\ 
%\vdots & & & \vdots \\
%\ol{\beta_{{n-1} 0}} & \ol{\beta_{{n-1} 1}} & \dots & \ol{\beta_{{n-1} n}}\\
%-1 & \alpha_0 &  \dots & \alpha_{n-1}
%\end{matrix} \right|
%= \ds \frac{-1}{\ol{D^{(n-1)}}} \left|
= \ds \frac{1}{\, {D^{(n-1)}}\,} \left|
\begin{matrix}
\beta_{0 \,0} & \beta_{1\, 0} & \dots & \beta_{n\, 0} \\ 
\vdots & & & \vdots \\
\beta_{0 \,{n-1}} & \beta_{1\, {n-1}} & \dots & \beta_{n \,{n-1}}\\
-1 & \alpha_0 &  \dots & \alpha_{n-1}
\end{matrix} \right|
\label{form1}
\ee
Let $C$ be the matrix with entries as in the determinant in (\ref{form1}) above. It could be expressed as follows
\be
C = \begin{pmatrix}
A & v \\ w & \alpha_{n-1}
\end{pmatrix}
\ee where $A$ is the $n \times n$ matrix with entries $A_{j k} = \beta_{k j}$, $v$ is the column vector $(\beta_{n 0}, \dots, \beta_{n {n-1}})^{T}$ and $w$ is the row vector $(-1, \alpha_0, \dots, \alpha_{n-2})$. Note that $\det(A)= {D^{(n-1)}}$ and it is real as $A$ is Hermitian.

Now we use Lemma \ref{lemma2} to compute $det(C)$. To do that, we need to find out what $A^{-1}$ is.

By the definition of $\nu$,
\be
A_{j k} = (1-\gamma) \|\Phi_k\|^2 \delta_{k j} +  \gamma \ol{\Phi_k(\zeta)} \Phi_j (\zeta) = \|\Phi_k\| \|\Phi_j\| M_{j k}
\ee where
\be
M_{j k} = (1-\gamma) \delta_{k j} +  \gamma \ol{\vp_k(\zeta)} \vp_j (\zeta)
\ee Observe that for any column vector $x =(x_0, x_1, \dots, x_{n-1})^{T}$,
\be
M x = (1-\gamma) x +  \gamma \left( \ds \sum_{j=0}^{n-1} \vp_j (\zeta) x_j \right) \left( \vp_0(\zeta), \vp_1(\zeta), \dots, \vp_0(\zeta) \right)^{T} 
\ee Therefore, if $P_{\vp}$ denotes the orthogonal projection onto the space spanned by the vector $\vp = \left( \vp_0(\zeta), \vp_1(\zeta), \dots, \vp_0(\zeta) \right)$, we can write
\be
M = (1-\gamma) \mathbf{1} + \gamma K_{n-1} P_{\vp}
\ee Hence, the inverse of $M$ is
\be
M^{-1} = (1-\gamma)^{-1} (\mathbf{1} - P_{\vp}) + ((1-\gamma)+ \gamma K_{n-1})^{-1} P_{\vp}
\ee and the inverse of $A$ is
\be
A^{-1} = D^{-1} M^{-1} D^{-1}
\ee where $D_{ij} = \|\Phi_i\| \delta_{i j}$.

Recall that $v = (\beta_{n 0}, \beta_{n 1}, \dots, \beta_{n {n-1}})^{T}$, which is a multiple of $\vp$. Therefore,
\be
(A^{-1} v)_j = \left( (1-\gamma) + \gamma K_{n-1} \right)^{-1} \gamma \hspace{0.05in} \ol{\Phi_n(\zeta)} \, \|\Phi_j\|^{-1} \vp_j(\zeta)
\label{form2}
\ee

(\ref{form2}), (\ref{form1}) and Lemma \ref{lemma2} then imply
\be
\alpha_{n-1}(d\nu) = \alpha_{n-1} - \left( (1-\gamma) + \gamma K_{n-1} \right)^{-1} \gamma \hspace{0.05in}  \ol{\vp_n(\zeta)} \left(\ds \sum_{j=0}^{n-1} \alpha_{j-1} \ds \frac{\|\Phi_n\|}{\|\Phi_j\|} \vp_j(z_0) \right)
\ee \end{proof}
This concludes the proof of Theorem \ref{barrytheorem}. \\

Now we are going to prove Theorem \ref{theorem0}.
\begin{proof}
First, observe that $\alpha_{j-1}=-\ol{\Phi_j(0)}$. Therefore, $\alpha_{j-1}/\|\Phi_j\| = - \ol{\vp_j(0)}$. Second, observe that $\|\Phi_{n+1}\|$ is independent of $j$ so it could be taken out from the summation. As a result, (\ref{mainformula}) in Theorem \ref{barrytheorem} becomes
\be
\alpha_{n}(d\nu)=\alpha_{n}(d\mu) + q_{n}^{-1} \gamma \hspace{0.05in} \ol{\vp_{n+1}(\zeta)}  \,\|\Phi_{n+1}\| \left( \ds \sum_{j=0}^{n} \ol{\vp_j(0)}\vp_j(\zeta) \right)
\label{intermediate}
\ee

Then we use the Christoffel-Darboux formula, which states that for $x, y \in \mathbb{C}$ with $x \bar{y} \not = 1$,
\begin{equation}
(1-\ol{x} y) \left( \ds \sum_{j=0}^{n} \ol{\vp_j(x)} \vp_j(y) \right) = \ol{\vp_{n}^*(x)}\vp_{n}^*(y) -  \ol{x}y \ol{\vp_n(x)} \vp_n (y)
\end{equation}

Moreover, note that $q_n^{-1} \gamma = ((1-\gamma)\gamma^{-1} + K_n(\zeta))^{-1}$
Therefore, (\ref{intermediate}) could be simplified as follows
\be 
\alpha_n(d\nu)
= \alpha_n + \ds \frac{\ol{\vp_{n+1}(\zeta)} \vp_n^*(0) \vp_n^*(\zeta)}{(1-\gamma)\gamma^{-1} + K_n(\zeta)} \|\Phi_{n+1}\|
\label{formula1}
\ee

Finally, observe that $\vp_n^*(0)=\|\Phi_n\|^{-1}$ and that by (\ref{normphi}), $\|\Phi_{n+1}\|/\|\Phi_n\|=(1-|\alpha_n|^2)^{1/2}$. This completes the proof.
\end{proof}

\section{Acknowledgements}
It is a pleasure to thank Professor Barry Simon for his guidance and Cherie Galvez for proof-reading this paper.


\bigskip

\begin{thebibliography}{8}

\bibitem{cm} A. Cachafeiro and F. Marcell\'an, \emph{Modifications of Toeplitz Matrices: jump functions}, Rocky Mountain J. Math. \textbf{23} (1993), 521-531

\bibitem{geronimus} Ya. L. Geronimus, \emph{Polynomials Orthogonal on a Circle and Their Applications}, Amer. Math. Soc. Translation \textbf{1954} (1954), no. 104, 79pp

\bibitem{nevai} P. Nevai, \emph{Orthogonal Polynomials}, Mem. Amer. Math. Soc. \textbf{18} (1979), no. 213, 185pp 

\bibitem{simon1} B. Simon, \emph{Orthogonal Polynomials on the Unit Circle, Part 1: Classical Theory}, AMS Colloquium Series, American Mathematical Society, Providence, RI, 2005.

\bibitem{simon2} B. Simon, \emph{Orthogonal Polynomials on the Unit Circle, Part 2: Spectral Theory}, AMS Colloquium Series, American Mathematical Society, Providence, RI, 2005.

%1026
\bibitem{szego} G. Szeg\H o, \emph{Orthogonal Polynomials}, Amer.\ Math.\ Soc.\ Colloq. Publ., Vol. 23, American Mathematical Society, Providence, RI, 1939; third edition, 1967.

\bibitem{wong1} M.-W. L. Wong, \emph{Generalized bounded variations and inserting point masses}, to appear in Constr. Approx.

\bibitem{wong2} M.-W. L. Wong, \emph{Asymptotics of polynomials and point perturbation in a gap}, preprint.


\end{thebibliography}
\end{document}